\documentclass[11pt,a4paper,reqno]{amsart} 
\usepackage{amsfonts,amsthm, amscd, epsfig, amsmath, amssymb,enumerate}
\newcommand{\halmos}{\rule{1ex}{1.4ex}}

\newcounter{for}[section]

\newtheorem{itlemma}{Lemma}[section]
\newtheorem{itproposition}[itlemma]{Proposition}
\newtheorem{itfact}[itlemma]{Fact}
\newtheorem{theorem}[itlemma]{Theorem}
\newtheorem{itcorollary}[itlemma]{Corollary}
\newtheorem{itremark}[itlemma]{Remark}
\newtheorem{itremarks}[itlemma]{Remarks}
\newtheorem{itdefinition}[itlemma]{Definition}
\newtheorem{itexample}[itlemma]{Example}

\newenvironment{fact}{\begin{itfact}\rm}{\end{itfact}}
\newenvironment{claim}{\begin{itclaim}\rm}{\end{itclaim}}
\newenvironment{lemma}{\begin{itlemma}}{\end{itlemma}}
\newenvironment{remark}{\begin{itremark}\rm}{\end{itremark}}
\newenvironment{remarks}{\begin{itremarks} \rm}{\end{itremarks}}
\newenvironment{corollary}{\begin{itcorollary}}{\end{itcorollary}}
\newenvironment{proposition}{\begin{itproposition}}{\end{itproposition}}
\newenvironment{definition}{\begin{itdefinition}\rm}{\end{itdefinition}}
\newenvironment{example}{\begin{itexample}\rm}{\end{itexample}}
\newcommand{\be}[1]{\addtocounter{for}{1} \begin{equation}\label{#1}}
\newcommand{\eeq}{\end{equation}}
\newcommand{\bl}[1]{\begin{lemma}\label{#1}}
\newcommand{\br}[1]{\begin{remark}\label{#1}}
\newcommand{\brs}[1]{\begin{remarks}\label{#1}}
\newcommand{\bt}[1]{\begin{theorem}\label{#1}}
\newcommand{\bd}[1]{\begin{definition}\label{#1}}
\newcommand{\bp}[1]{\begin{proposition}\label{#1}}
\newcommand{\bfact}[1]{\begin{fact}\label{#1}}
\newcommand{\bc}[1]{\begin{corollary}\label{#1}}
\newcommand{\bex}[1]{\begin{example}\label{#1}}
\newcommand{\ec}{\end{corollary}}
\newcommand{\efact}{\end{fact}}
\newcommand{\eex}{\end{example}}
\newcommand{\el}{\end{lemma}}
\newcommand{\er}{\end{remark}}
\newcommand{\ers}{\end{remarks}}
\newcommand{\et}{\end{theorem}}
\newcommand{\ed}{\end{definition}}
\newcommand{\ep}{\end{proposition}}
\newcommand{\epr}{\end{proof}}
\newcommand{\bpr}{\begin{proof}}
\newcommand{\bcl}[1]{\begin{claim}\label{#1}}
\newcommand{\ecl}{\end{claim}}

\newcommand{\ecs}{\end{corollary}}
\newcommand{\eers}{\end{exercise}}
\newcommand{\eexs}{\end{example}}
\newcommand{\eems}{\end{example}}
\newcommand{\els}{\end{lemma}}
\newcommand{\eles}{\end{lemmaex}}
\newcommand{\ets}{\end{theorem}}
\newcommand{\eds}{\end{definition}}
\newcommand{\eps}{\end{proposition}}

\newcommand{\bi}{\begin{itemize}}
\newcommand{\ei}{\end{itemize}}
\newcommand{\ben}{\begin{enumerate}}
\newcommand{\een}{\end{enumerate}}

\def\vbar{\mathchoice{\vrule height6.3ptdepth-.5ptwidth.8pt\kern-.8pt}
   {\vrule height6.3ptdepth-.5ptwidth.8pt\kern-.8pt}
   {\vrule height4.1ptdepth-.35ptwidth.6pt\kern-.6pt}
   {\vrule height3.1ptdepth-.25ptwidth.5pt\kern-.5pt}}
\def\fudge{\mathchoice{}{}{\mkern.5mu}{\mkern.8mu}}
\def\bbc#1#2{{\rm \mkern#2mu\vbar\mkern-#2mu#1}}
\def\bbb#1{{\rm I\mkern-3.5mu #1}}
\def\bba#1#2{{\rm #1\mkern-#2mu\fudge #1}}
\def\bb#1{{\count4=`#1 \advance\count4by-64 \ifcase\count4\or\bba A{11.5}\or
   \bbb B\or\bbc C{5}\or\bbb D\or\bbb E\or\bbb F \or\bbc G{5}\or\bbb H\or
   \bbb I\or\bbc J{3}\or\bbb K\or\bbb L \or\bbb M\or\bbb N\or\bbc O{5} \or
   \bbb P\or\bbc Q{5}\or\bbb R\or\bbc S{4.2}\or\bba T{10.5}\or\bbc U{5}\or
   \bba V{12}\or\bba W{16.5}\or\bba X{11}\or\bba Y{11.7}\or\bba Z{7.5}\fi}}

\def \qed {{\hspace*{\fill}$\halmos$\medskip}}

\def \Z {{\mathbb Z}}
\def \R {{\mathbb R}}

\def \ra {\rightarrow }

\def \o {\omega}

\def \O{\Omega}
\def \s{\sigma}

\def\g{\gamma}
\def\d{\delta}
\def\e{\varepsilon}
\def\S{\Sigma}
\def\b{\beta}
\def\L{\Lambda}
\def\l{\lambda}
\def\E{{\mathcal{E}}}
\def\D{{\mathcal{D}}}
\def\1{{\bf 1}}

\numberwithin{equation}{section}

\newtheorem{teo}{Theorem}[section]

\newtheorem{de}[teo]{Definition}

\DeclareMathOperator{\distanza}{dist}

\renewcommand{\tilde}{\widetilde}

\newcommand{\ie}{\emph{i.e.\ }}

\DeclareMathSymbol{\leqslant}{\mathalpha}{AMSa}{"36} 
\DeclareMathSymbol{\geqslant}{\mathalpha}{AMSa}{"3E} 
\DeclareMathSymbol{\eset}{\mathalpha}{AMSb}{"3F}     
\renewcommand{\leq}{\;\leqslant\;}                   
\renewcommand{\geq}{\;\geqslant\;}                   

\newcommand{\grad}{\nabla}

\def\1{\ifmmode {1\hskip -3pt \rm{I}} \else {\hbox {$1\hskip -3pt \rm{I}$}}\fi}

 \let\b=\beta  \let\c=\chi \let\d=\delta  \let\e=\varepsilon
 \let\g=\gamma \let\h=\eta      \let\l=\lambda
      \let\o=\omega      
  \let\s=\sigma \let\t=\tau   
  
\let\D=\Delta     \let\L=\Lambda 
\let\O=\Omega      

\def\\{\hfill\break}

\def\tthsp{\kern .083333 em}

\def\?{\mskip -10mu}
\def\indbox#1{\hbox to \parindent{\hfil\ #1\hfil} }

\newfam\msafam
\newfam\msbfam
\newfam\eufmfam
\def\hexnumber#1{%
  \ifcase#1 0\or 1\or 2\or 3\or 4\or 5\or 6\or 7\or 8\or
  9\or A\or B\or C\or D\or E\or F\fi}
\font\tenmsa=msam10
\font\sevenmsa=msam7
\font\fivemsa=msam5
\textfont\msafam=\tenmsa
\scriptfont\msafam=\sevenmsa
\scriptscriptfont\msafam=\fivemsa        
\edef\msafamhexnumber{\hexnumber\msafam}%
\mathchardef\restriction"1\msafamhexnumber16
\mathchardef\ssim"0218
\mathchardef\square"0\msafamhexnumber03
\mathchardef\eqd"3\msafamhexnumber2C
\def\QED{\ifhmode\unskip\nobreak\fi\quad
  \ifmmode\square\else$\square$\fi}            
\font\tenmsb=msbm10
\font\sevenmsb=msbm7
\font\fivemsb=msbm5
\textfont\msbfam=\tenmsb
\scriptfont\msbfam=\sevenmsb
\scriptscriptfont\msbfam=\fivemsb
\def\Bbb#1{\fam\msbfam\relax#1}    
\font\teneufm=eufm10
\font\seveneufm=eufm7
\font\fiveeufm=eufm5
\textfont\eufmfam=\teneufm
\scriptfont\eufmfam=\seveneufm
\scriptscriptfont\eufmfam=\fiveeufm

\def\A{{\mathcal{A}}}

\def\({\left(}
\def\){\right)}

\def\o#1{{\rm [\mkern -3mu #1}}
\def\c#1{{\rm ]\mkern -3mu #1}}
\let\Z=\integer

\let\neper=e
\let\ii=i

\def\ie{\hbox{\it i.e.\ }}

\let\<=\langle
\let\>=\rangle

\def\Tr{ \mathop{\rm Tr}\nolimits }

\outer\def\nproclaim#1 [#2]#3. #4\par{\medbreak \noindent
   \talato(#2){\bf #1 \Thm[#2]#3.\enspace }%
   {\sl #4\par }\ifdim \lastskip <\medskipamount 
   \removelastskip \penalty 55\medskip \fi}
\def\thmm[#1]{#1}
\def\teo[#1]{#1}
\def\sttilde#1{%
\dimen2=\fontdimen5\textfont0
\setbox0=\hbox{$\mathchar"7E$}
\setbox1=\hbox{$\scriptstyle #1$}
\dimen0=\wd0
\dimen1=\wd1
\advance\dimen1 by -\dimen0
\divide\dimen1 by 2
\vbox{\offinterlineskip%
   \moveright\dimen1 \box0 \kern - \dimen2\box1}
}
\newcommand{\Hamcl}{\mathcal{H}_{\mbox{cl}}}
\newcommand{\Ham}{\mathcal{H}}

\newcommand{\ket}[1]{|\,#1 \,\rangle}

\newcommand{\braket}[3]{\langle #1\,|\,#2\,|\,#3 \rangle}

\newcommand{\De}{\mathrm{d}}

\newcommand{\mnorm}[1]{%
   \left\vert\mkern-1mu\left\vert\mkern-1mu\left\vert #1
     \right\vert\mkern-1mu\right\vert\mkern-1mu\right\vert}

\newcommand{\gen}{\mathcal{L} }
\newcommand{\genh}{\hat{\mathcal{L} }}

\newcommand{\expo}[1]{\exp{\left[#1\right]}}

\newcommand{\hS}{\hat{S}}
\newcommand{\hp}{\hat{\pi}}
\newcommand{\hP}{\hat{\mathbb{P}}}
\newcommand{\Hb}{\hat{H}_{\beta}}

\begin{document}

\title[Decay of correlations]
{Decay of correlations for quantum spin systems with a transverse field: a dynamic approach}

\begin{abstract}
We consider a wide class of quantum spin systems obtained by adding a transverse field to a classical Hamiltonian. We give explicit high-temperature conditions which guarantee exponential decay of correlations. A stochastic-geometric representation allows to reformulate the model as a classical random field on a space of marked point processes on the circle $[0,\b)$, where $\b$ is the inverse temperature. We then construct a Markov process having this random field as invariant measure. By the mixing properties of the process, the exponential decay of correlations follows by an adaptation of a general argument.



\end{abstract}
\author[A. Cipriani]{Alessandra Cipriani}
\address{Institut f\"{u}r Mathematik,  Universit\"{a}t Z\"{u}rich,
Winterthurerstrasse 190, CH-8057,
Z\"{u}rich, Switzerland} \email{alessandra.cipriani@math.uzh.ch}
\author[P. Dai Pra]{Paolo Dai Pra}
\address{Dipartimento di Matematica Pura e Applicata, 
Universit\`{a} di Padova, Via Belzoni 7, 35131 Padova, Italy}
\email{daipra@math.unipd.it}



\maketitle

\thispagestyle{empty}

\section{Introduction}

It is well known that many spatial properties of random fields can be
detected by studying properties of some Markov process which is
reversible with respect to the distribution of the field. For example,
as established for many lattice spin systems (\cite{ma:ol,Yos}) and for particle Glauber dynamics in the continuum (\cite{BCC}), mixing properties of the
field correspond to fast convergence to equilibrium of the associated
Markovian dynamics. For many purposes, the right notion of fast
convergence to equilibrium is equivalent to the fact that the generator
$\gen$ of the dynamics has a strictly positive {\em spectral gap}. Thus,
whenever sufficient conditions for positivity of the spectral gap are
available, one obtains that the field is spatially mixing in a rather
strong sense. 

The aim of this paper is to show that this approach can be applied to
a class of {\em quantum fields} obtaining, as final consequence,
sufficient conditions for exponential decay of quantum correlations. We
consider Gibbs fields on a lattice $\Lambda$ obtained from a classical
Hamiltonian by adding a {\em transverse field}. Our approach consists of
the following steps.

\begin{enumerate}
\item
We use a stochastic-geometric representation of this fields that allows
to reformulate quantum means and correlations of observables as {\em
classical} means and correlations with respect to a suitable
distribution $\pi$ on a product space of {\em marked point processes}. This representation dates back to \cite{AizNacht}, and it has been recently developed and applied in \cite{Ioffe, CCIL}.
\item
We then build a Markovian dynamics that is reversible with respect to
$\pi$ (actually, to a suitable {\em lifting} of $\pi$). In this context, the Bochner-Bakry-Emery
method (or $\Gamma_2$ method) for estimating the spectral gap, in the
formulation for discrete settings given in \cite{BCDP}, happens to be
quite successful. In particular, we obtain simple and explicit high
temperature conditions for positivity of the spectral gap.
\item
By a nontrivial adaptation to this setting of a classical argument (see
e.g. \cite{ma, BCC}), we show that positivity of the spectral gap implies
exponential decay of correlations under $\pi$, which by the stochastic-geometric representation
implies decays of correlations for the original quantum model. 
\end{enumerate}

The most used method for estimating correlations in high temperature regimes, for both classical and quantum models, is {\em cluster expansion} (see e.g. \cite{PoZe, Uel} for some rigorous results in this direction). Therefore, the approach proposed in this paper should be compared with cluster expansion. Our method is very indirect, it requires the introduction of stochastic dynamics and rather sophisticated abstract tools to study them. On the other hand, cluster expansion is a very direct combinatorial method. However, when one tries to obtain explicit estimates on a specific class of models, our method has good chances to lead to simple expressions, so that explicit high-temperature conditions can often be obtained. Moreover, what we propose is quite robust in terms of the structure of the Hamiltonian of the systems; in particular it works well in presence of spatial inhomogeneities.

The paper is organized as follows. In Section 2 we define the class of models we consider, and state our result on decay of correlations. In the following two sections we review the stochastic geometric representation in terms of marked point processes (Section 3) and the discrete Bochner-Bakry-Emery
method for estimating the spectral gap of a Markov generator (Section 4). In Section 5 we introduce the Markovian dynamics which corresponds to the representation in Section 3, and we prove that the corresponding generator has a uniformly (in the volume) positive spectral gap. Finally, in Section 6 we show that uniform positivity of the spectral gap implies exponential decay of correlations.

\section{Model and main result}

In correspondence with the classical spin values $\pm 1$, we consider the two dimensional vectors
\[
\ket{1}:={1 \choose 0} \qquad \ket{-1}:={0 \choose 1}.
\]
The pair $\ket{1}$ and $\ket{-1}$ is a basis of the state space for a one-site quantum spin system, provided with the Euclidean scalar product $\langle\,\cdot \, | \, \cdot \,\rangle_2$.
For a finite $\L \subseteq \Z^d$, we consider the classical  space $\O_{\L} := \{-1,1\}^{\L}$. The corresponding quantum space is given by
\[
\mathbb{X}_\Lambda :=\bigotimes_{i \in \Lambda} \mathbb{C}^2.
\]
We denote by $\ket{\sigma}:= \otimes_{i \in \Lambda} \ket{\sigma_i}$ the $2^{|\L|}$ elements of a basis of $\mathbb{X}_\Lambda$; the scalar product, defined on the elements of the basis by
\[
\langle\,\sigma\, | \, \sigma'\,\rangle := \prod_{i \in \Lambda}^{}\langle\,\sigma_i\, | \, {\sigma'}_i\,\rangle_2 ,
\]
makes it into an orthonormal basis.

Let now by $\hat{\sigma}^z$ the Pauli matrix, \[
\hat{\sigma}^z:= \left(\begin{array}{cc} 1 & 0 \\ 0 & -1 \end{array}\right),
\]
and denote by $\hat{\sigma}^z_i$ the linear operator on $\mathbb{X}_\Lambda$ defined on the basis by
\be{pauliz}
\hat{\sigma}^z_i\ket{\sigma}:= \ket{\sigma_1} \otimes \ldots \otimes  \hat{\sigma}^z \ket{\sigma_i}\otimes \ldots \ket{\sigma_{|\Lambda|}}.
\eeq
Given a function $H: \O_{\L} \ra \R$, it can be lifted to the self-adjoint operator $\Hamcl$ on $\mathbb{X}_\Lambda$ by
\[
\Hamcl \ket{\sigma} = H(\hat{\sigma}^z_1, \ldots, \hat{\sigma}^z_{|\L|}) \ket{\sigma} = H(\s) \ket{\sigma} .
\]
The models we consider are obtained by adding a {\em transverse field} to a given classical Hamiltonian $\Hamcl$. To the Pauli matrix
$$\hat{\sigma}^x:= \left(\begin{array}{cc} 0 & 1\\ 1 & 0 \end{array} \right)$$
we associate, as in \eqref{pauliz}, the operators $\hat{\sigma}^x_i$. Given real numbers $\{\l_i: \, i \in \L \}$, we define the Hamiltonian
\be{hamiltonian}
\Ham := \Hamcl - \sum_{i \in \L} \l_i \hat{\sigma}^x_i.
\eeq
\br{positivity}
It should be noticed that there is some arbitrariness in the choice of the basis $\ket{1}, \ket{-1}$ of eigenvectors for $\hat{\s}^z$. For instance, if we transform the canonical choice given above to $\ket{1}:={-1 \choose 0}$ and $ \ket{-1}:={0 \choose 1}$, the corresponding change of variables would transform $\hat{\s}^x$ to $-\hat{\s}^x$. Since this choice can be made componentwise, we can always choose a basis of $\mathbb{X}_\Lambda$ such that the Hamiltonian corresponds to the matrix (\ref{hamiltonian}) {\em with nonnegative values of the} $\l_i$'s. Thus, by no loss of generality, we will assume in the rest of the paper $\l_i \geq 0$ for every $i \in \L$.
\er

We recall that for an observable $\mathcal F$ its {\em quantum average} is given by
$$<\mathcal{F}>:= \frac{\Tr{\mathcal{F} e^{-\beta \Ham} }}{\Tr{ e^{-\beta \Ham} }}. $$
Similarly, for two observables $\mathcal F$ and $\mathcal G$, the {\em (truncated) correlation} is defined by
$$<\mathcal F; \mathcal G>:=<\mathcal F \mathcal G>-<\mathcal F>< \mathcal G>.$$
We will be interested, in particular, in observations of the spin. In other words, to the classical observable $f:\O_{\L} \ra \R$, we associate the quantum counterpart
\[
\mathcal{F}_f := f(\hat{\sigma}^z_1, \ldots, \hat{\sigma}^z_{|\L|}).
\]
Before stating our main result, we need to introduce some further notations and assumptions. For $\s \in \O_{\L}$ and $i \in \L$, $\s^i$ will denote the element of $\O_{\L}$ obtained from $\s$ by flipping the $i^{\mbox{th}}$ spin; also, for $f :\O_{\L} \ra \R$, we set
\[
\grad_i f(\s) := f(\s^i) - f(\s).
\]
Moreover, we define
\be{norms}
\|f\| := \sup_{\s \in \O_{\L}} |f(\s)| \hspace{1cm} \mnorm{f} := \sum_{i \in \L} \|\grad_i f\|.
\eeq
Finally, we let
\be{support}
\L_f := \{ i \in \L: \|\grad_i f\| >0\}.
\eeq
In other words, $i \not\in \L_f$ if $f$ does not depend on the spin $\s_i$.

From the Statistical Mechanics point of view, one is interested in the infinite volume limit $\L \uparrow \Z^d$. Thus $\L$ should be interpreted as a parameter in our model in finite volume: it will be relevant to specify, both in the assumptions and in the results, which constants do {\em not} depend on $\L$, and are therefore preserved in the infinite volume limit, whenever such limit makes sense. In this paper we assume that the classical part of the Hamiltonian satisfies the following properties, for some constants $C,\l,R>0$ that do not depend on $\L$.
\bi
\item[({\bf BD})]
\[
\sup_{i \in \L} \sup_{\s \in \O_{\L}} \left| \grad_i H(\s) \right| \leq C \hspace{1.5cm} 0 \leq  \l_i \leq \l \ \ \forall i \in \L.
\]
\item[({\bf LOC})]
For every $i \in \L$
\[
\grad_i H(\s) = \grad_i H(\eta)
\]
whenever $\s_j = \eta_j$ for every $j$ such that $|j-i| \leq R$.
\ei

We now state the main result of this paper.
\bt{decay}
Assume  properties ({\bf BD}) and ({\bf LOC}) hold for $H$, and suppose $\b,\l$ are such that
\[
\g := e^{-2\b C}  -  5 \max(1,\b \l) (2R)^d \left(e^{3 \b C} -1\right) >0
\]
(note that, for $\l$ given, this is true for $\b$ sufficiently small). Then
\[
\left|<{\mathcal F}_f; {\mathcal F}_g> \right| \leq \left(\|f\| \|g\| + N \mnorm{f}\mnorm{g} \right) \exp[- \d \distanza(\L_f,\L_g)],
\]
where $\distanza(\L_f,\L_g) := \min\{|i-j| : i \in \L_f, \, j \in \L_g\}$, and
\[
N = \left(\frac{2}{1-e^{-1/R}}\right)^d, \hspace{2cm} \d := \frac{1}{2R} \min\left(\frac{1}{2}, \frac{\g}{10 \max(1,\l \b) e^{\b C+1}}\right).
\]
\et

\section{Stochastic-geometric representations}

We begin by introducing some notations. We denote by $\D$ the set of piecewise constant, right continuous functions from $[0,\b) \ra \{-1,1\}$, where $[0,\b)$ is meant to be a circle ($0=\b$). This set is provided with the Skorohod topology. Moreover, $\S$ is the set of finite subsets of $[0,\b)$. The topology on $\S$ is generated by the following sets, parametrized by $\eta \in \S$ and $ \e>0$: $\{\eta' \in \S : |\eta'| = |\eta| \mbox{ and } \distanza(\eta,\eta') < \e\}$, where $\distanza(\eta,\eta')$ denotes the Hausdorff distance between two sets.
Finally, we define
\[
S := (\D \times \S)^{\L},
\] 
and provide it with the product topology.

If $\eta \in \S$ and $a \in \D$, we say that $a \sim \eta$ ($a$ is compatible with $\eta$) if the discontinuity points of $a$ are a subset of $\eta$. Similarly, for $\xi \in \S^{\L}$ and $\s \in \D^{\L}$ we write, by slight abuse of notation, $\s \sim \xi$ if $\s_i \sim \xi_i$ for every $i \in \L$. We will sometimes say that $\s$ is a {\em coloring} of $\xi$.

By $P_i$, $i \in \L$, we denote the Poisson point measure on $\S$ with intensity $\l_i$. By $\mathbb P$ we mean the product measure on $\S^{\L}$
\[
\mathbb P = \otimes_{i \in \L} P_i.
\] 
Consider the probability measure on $S$ given by
\be{pi}
\pi(\s,d\xi) := \frac{1}{Z} {\bf 1}_{\{\s \sim \xi\}} \exp\left[ \int_0^{\b} H(\s(t))dt \right] \mathbb P(d\xi),
\eeq
where $Z$ is a normalization factor.
Given Borel measurable functions $\Phi,\Psi: S \ra \R$, we denote by
\[
\pi[\Phi] \ \ \mbox{ or } \ \ \pi[\Phi(\s,\xi)]
\]
the mean of $\Phi$ with respect to $\pi$, and by
\[
\pi[\Phi;\Psi] \ \ \mbox{ or } \ \ \pi[\Phi(\s,\xi);\Psi(\s,\xi)]
\]
their covariance $\pi[\Phi \Psi] - \pi[\Phi]\pi[\Psi]$.

The following Theorem is a special case of what shown in \cite{Ioffe}, Section 1.2. The proof is sketched here for completeness.

\bt{stoch-geom}
Let $f,g : \O_{\L} \ra \R$. Then
\be{reprmean}
<{\mathcal F}_f> = \pi\left[f(\s(0))\right],
\eeq
and
\be{reprcov}
<{\mathcal F}_f; {\mathcal F}_g> = \pi \left[f(\s(0));g(\s(0)) \right].
\eeq
\et
\bpr
We begin by giving a representation of the partition function $\Tr{ e^{-\beta \Ham} }$. First we set $K_i := \hat{\sigma}^x_i + \mathbb{I}$, and write
\[
e^{- \beta \Ham} =  \exp\left[-\beta \Hamcl + \beta\sum_i \lambda_i \hat{\sigma}^x_i \right] 
= \expo{-\beta \Hamcl+\beta \sum_i \left(\lambda_i K_i-\lambda_i\mathbb{I}\right)} 
\]
By the {\em Lie-Trotter formula}
\begin{multline*}
e^{- \beta \Ham}= \lim_{N \ra +\infty} \left(\expo{- \frac{\b}{N} \Hamcl}  \prod_{i \in \L} \expo{\frac{\b}{N} \l_i (K_i - \mathbb{I})} \right)^N \\ = \lim_{N \ra +\infty} \left(\expo{- \frac{\b}{N} \Hamcl}  \prod_{i \in \L} \left[ \left(1-\frac{\b \l_i}{N} \right) \mathbb{I} + \frac{\b \l_i}{N} K_i \right] \right)^N.
\end{multline*}
Now, let $X_i(k)$, $i \in \L, k \in \{1,2,\ldots,N\}$, be independent, $\{0,1\}$-valued random variables, such that $P(X_i(k) =1) = \frac{\b \l_i}{N}$. If ${\bf E}_N(\cdot)$ denotes expectation with respect to the joint law of the $X_i(k)$'s, we can write, for each $k$
\[
 \prod_{i \in \L} \left[ \left(1-\frac{\b \l_i}{N} \right) \mathbb{I} + \frac{\b \l_i}{N} K_i \right] = {\bf E}_N \left[\prod_{i \in \L} [(1-X_i(k)) \mathbb{I} + X_i(k) K_i] \right],
 \]
where the expectation of a matrix-valued random variable is just defined componentwise. Taking the product over $k$, we have
\begin{multline}
 \left(\expo{- \frac{\b}{N} \Hamcl}  \prod_{i \in \L} \expo{\frac{\b}{N} \l_i (K_i - \mathbb{I})} \right)^N  \\ =
 {\bf E}_N \left[ \prod_{k=1}^N \expo{- \frac{\b}{N} \Hamcl} 
\prod_{i \in \L} [(1-X_i(k)) \mathbb{I} + X_i(k) K_i] \right].
 \label{prod1}
 \end{multline}
To make the previous formula more compact, define $A_k(X) := \expo{- \frac{\b}{N} \Hamcl} 
\prod_{i \in \L} [(1-X_i(k)) \mathbb{I} + X_i(k) K_i]$. Note that
\[
A_k(X)  =  
\expo{- \frac{\b}{N} \Hamcl} \prod_{i:X_i(k) = 1} K_i
\]
We then obtain
\be{prod2}
\Tr{\left(\expo{- \frac{\b}{N} \Hamcl}  \prod_{i \in \L} \expo{\frac{\b}{N} \l_i (K_i - \mathbb{I})} \right)^N} = {\bf E}_N \left[\Tr{ \prod_{k=1}^N A_k(X)} \right].
\eeq
On the other hand
\[
\Tr{ \prod_{k=1}^N A_k(X)}  
= \sum_{\eta \in \O_{\L}} \braket{\eta}{\prod_{k=1}^N A_k(X)}{\eta} 
= \sum_{\eta(0),\eta(1),\ldots,\eta(N-1)} \prod_{k=1}^N \braket{\eta(k-1)}{A_k(X)}{\eta(k)},
\]
where we mean $\eta(N) = \eta(0)$. By definition of $K_i$, it is not hard to check that
\[
\braket{\eta(k-1)}{A_k(X)}{\eta(k)} = 
\begin{cases} 
\expo{- \frac{\b}{N}H(\eta(k-1))} & \mbox{if }  \eta_i(k-1) = \eta_i(k)  \mbox{ whenever } X_i(k) = 0 \\ 0 & \mbox{otherwise}.
\end{cases}
\]
Now we reformulate the above expressions in different terms. We identify the family of random variables $X=(X_i(k))$ with the point process $\xi^{(N)} = (\xi^{(N)}_i)_{i \in \L}$ on $[0,\b)^{\L}$ 	 by
\[
t \in \xi^{(N)}_i \ \iff \ t = \b \frac{k}{N} \mbox{ and } 	X_i(k) =1.
\]
Moreover, recalling what defined at the beginning of this section, $\s \in \D^{\L}$ is said to be {\em compatible} with $\xi^{(N)}$ ($\s \sim \xi^{(N)}$) if the discontinuities of $\s$ are a subset of $	\xi^{(N)}_i$. We obtain
\[
\sum_{\eta(0),\eta(1),\ldots,\eta(N-1)} \prod_{k=1}^N \braket{\eta(k-1)}{A_k(X)}{\eta(k)} = \sum_{\s \sim \xi^{(N)}} \expo{- \frac{\b}{N} \sum_{k=1}^N H(\s(\b k/N))}.
\]
Since each $\s \sim \xi^{(N)}$ is constant on each interval of the form $\left[ \frac{\b k}{N}, \frac{\b(k+1)}{N} \right)$, this last expression equals
\[
\sum_{\s \sim \xi^{(N)}} \expo{-  \int_0^{\b} H(\s(t))dt}.
\]
Summing up, we have shown that
\be{partfun1}
\Tr{ e^{-\beta \Ham} } = \lim_{N \ra +\infty} \mathbb{E} \left[ \sum_{\s \sim \xi^{(N)}} \expo{-  \int_0^{\b} H(\s(t))dt} \right],
\eeq
where the expectation is over the law of $\xi^{(N)}$. Now it is well known that the point process $\xi^{(N)}$ converges weakly to $\xi = (\xi_i)_{i \in \L}$ where the $\xi_i$'s are independent Poisson processes with intensities $\l_i$. The passage to the limit in (\ref{partfun1}) needs, however, to be justified, since
\[
F(\xi) :=  \sum_{\s \sim \xi} \expo{-  \int_0^{\b} H(\s(t))dt} 
\]
is continuous but not bounded in $\xi$. It is immediately seen, however, that $F(\xi) \leq C2^{|\xi|}$ for some $C>0$. Since, for the limiting Poisson process, the law of $|\xi|$ is Poissonian, and has therefore tails smaller than exponential, a standard truncation argument applies (we omit the details), and we get
\be{partfun}
\Tr{ e^{-\beta \Ham} } =\mathbb{E} \left[ \sum_{\s \sim \xi} \expo{-  \int_0^{\b} H(\s(t))dt} \right].
\eeq

This argument extends readily to traces of the form $\Tr{ \mathcal{F}_f e^{-\beta \Ham} }$, the only difference being that in the right hand side of (\ref{prod2}) we get
\[
{\bf E}_N \left[\Tr{  \mathcal{F}_f \prod_{k=1}^N A_k(X)} \right].
\]
The same product expansion leads to
\[
\Tr{ \mathcal{F}_f e^{-\beta \Ham} } = \mathbb{E} \left[ \sum_{\s \sim \xi} f(\s(0))\expo{-  \int_0^{\b} H(\s(t))dt} \right],
\]
which completes the proof of (\ref{reprmean}). Since, for the same reason, we also have
\[
\Tr{ \mathcal{F}_f \mathcal{F}_g e^{-\beta \Ham} } = \mathbb{E} \left[ \sum_{\s \sim \xi} f(\s(0)) g(\s(0))\expo{-  \int_0^{\b} H(\s(t))dt} \right],
\]
we also obtain (\ref{reprcov}).
\epr

\section{Spectral Gap of a Markov generator via the Bochner-Bakry-Emery method}

In this section we summarize the version of the $\Gamma_2$ method (see e.g. \cite{BakryEmery}) for the estimate of the spectral Gap of a Markov generator, as developed in \cite{BCDP} for processes with jumps. We consider Markov processes taking values in a measurable space $(S,\mathcal{S})$. Let $G$ be a subset of $S^S$, provided with a $\s$-field $\mathcal{G} \subseteq \mathcal{P}(G)$, having the property that the map $(\eta,\g) \mapsto \g(\eta)$, from $S \times G$ to $G$ is measurable (product spaces are meant to be provided with the product $\s$-field). We then associate to $\eta \in S$ a $\sigma$-finite positive measure (denoted by $c(\eta,\De \gamma)$) on $(G,\mathcal{G})$. This measure must satisfy the condition that, for every measurable $\phi:G \to[0,+ \infty]$, the function $\eta \mapsto \int \phi(\gamma)c(\eta,\De\gamma)$ is measurable as well. We assume that the considered Markov process has an infinitesimal generator $\gen$ that can be written in the following form:
\begin{equation}\label{eq:Lintegrale}
\mathcal{L}f(\eta)=\int_G c(\eta, \De \gamma)\nabla_\gamma f(\eta) ,
\end{equation}
for $f \in \mathcal{D}(\mathcal{L})$, having set $\nabla_\gamma f:=f \circ \gamma- f$. For simplicity, we assume that bounded measurable functions form a {\em core} for $\gen$.

  Let now $\pi$ be a probability measure on $(S,\mathcal{S})$, and denote by $\pi_c$  the positive measure on $S \times G$ given by $\pi_c(d\eta,d\g) := \pi(d\eta)c(\eta,d\g)$. We make the following additional assumptions on the generator ${\mathcal{L}}$.

\bi
\item[{\bf (Rev)}]
 For every $\g \in G$ there is a unique $\g^{-1} \in G$ such that the equality $\g^{-1}(\g(\eta)) = \eta$ holds $\pi_c$-a.s., and the map $\g \mapsto \g^{-1}$ is a measurable bijection on $G$. Moreover, for every $\Psi : S \times G \ra [0,+\infty)$ measurable, 
 \ei
\be{db}
\int \Psi(\eta,\g) c(\eta,d\g) \pi(d\eta) = \int \Psi(\g(\eta),\g^{-1}) c(\eta,d\g) \pi(d\eta).
\eeq

\vspace{0.5cm}
Note that assumption {\bf (Rev)} implies that ${\mathcal L}$ is symmetric in $L^2(\pi)$, \ie
\be{sym}
\int f(\eta)g(\g(\eta))c(\eta,d\g)\pi(d\eta) = \int f(\g(\eta))g(\eta)c(\eta,d\g)\pi(d\eta) 
\eeq
for $f,g \in \mathcal{D}(\mathcal{L})$. Then, being a symmetric Markov operator, $\mathcal{L}$ is self-adjoint in $L^2(\pi)$. Thus, {\bf (Rev)} is a reversibility condition,  (\ref{db}) is the usual {\em detailed balance} condition written in this general context and $\pi$ is a stationary distribution for the process. 
\\ Note that, under  {\bf (Rev)}, for $f,g \in {\mathcal D}({\mathcal L})$,
\be{dir}
\E(f,g) := - \pi\left(f{\mathcal L} g\right) = \frac{1}{2}\pi\left[ \int c(\eta,d\g) \nabla_{\g}f(\eta) \nabla_{\g}g(\eta)\right].
\eeq
$\E(f,g)$ is called the {\em Dirichlet form} associated to the process. 

We recall that the {\em spectral Gap}  gap$(\mathcal{L})$ of $\mathcal{L}$ is defined as the largest $k>0$ such that the Poincar\'e inequality
\[
k \pi[f;f] \leq \E(f,f)
\]
holds for every $f \in \mathcal{D}(\mathcal{L})$. An alternative representation of the spectral Gap is given in the following Proposition (see \cite{BCDP}, Proposition 1.1).
\bp{Gamma2}
The spectral gap gap$(\mathcal{L})$ equals the largest constant $k$ such that
\[
k \E(f,f) \leq \pi\left[ \left( \mathcal{L} f \right)^2 \right].
\]
\ep
The following definitions provides the key notion for the use of Proposition \ref{Gamma2} in the estimation of spectral gaps.
\bd{symmetric}
A measurable function $r: S \times G \times G \ra [0,+\infty)$ is said {\em $\gen$-symmetric} if the following conditions hold
\bi
\item[{\bf (A1)}]
$r \in L^1(\pi(d\eta) c(\eta,d\g) c(\eta,d\d))$.
\item[{\bf (A2)}]
The equality $r(\eta,\g,\d) = r(\eta,\d,\g)$ holds $\pi(d\eta) c(\eta,d\g) c(\eta,d\d)$-almost surely.
\item[{\bf (A3)}]
$r(\eta,\g,\d) >0 \ \Rightarrow \ \g(\d(\eta)) = \d(\g(\eta))$.
\item[{\bf (A4)}]
For every $F: S \times G \times G \ra [0,+\infty)$ measurable, we have
\[
\int F(\eta,\g,\d)r(\eta,\g,\d) \pi(d\eta) c(\eta,d\g) c(\eta,d\d) = \int F(\g(\eta),\g^{-1},\d)r(\eta,\g,\d) \pi(d\eta) c(\eta,d\g) c(\eta,d\d).
\]
\ei
\ed
 The following Theorem, proved in \cite{BCDP} (Corollaries 2.2 and 2.3) is the main tool we shall use in next section.
\bt{gap}
Let $r : S \times G \times G \ra [0,+\infty)$ be a $\gen$-symmetric function. Then, for every
$f:S \ra \R$ measurable and bounded, the inequality
\[
\pi\left[ \left( \mathcal{L} f \right)^2 \right] \geq \int \pi(d\eta)c(\eta,d\g)
c(\eta,d\d)[1-r(\eta,\g,\d)] \nabla_{\g}f(\eta)\nabla_{\d}f(\eta)
\]
holds. Therefore (see Proposition \ref{Gamma2}), if
\[
\int \pi(d\eta)c(\eta,d\g)
c(\eta,d\d)[1-r(\eta,\g,\d)] \nabla_{\g}f(\eta)\nabla_{\d}f(\eta) \geq k \E(f,f),
\]
then gap$(\mathcal{L}) \geq k$.
\et
The construction of an effective $\gen$-symmetric function is model dependent. We will construct one for the specific model in next section.

\section{The stochastic dynamics}

\subsection{Modification of the state space}

In (\ref{pi}) we have defined the probability measure
\[
\pi(\s,d\xi) := \frac{1}{Z} {\bf 1}_{\{\s \sim \xi\}} \exp\left[ \int_0^{\b} H(\s(t))dt \right] \mathbb P(d\xi)
\]
on $S := (\D \times \S)^{\L}$.  It will be convenient, for our purposes, to reformulate the model in a different, but equivalent, state space. We first note that $\pi$ is concentrated on the set $\{ (\s,\xi) \in S : \s \sim \xi\}$. Given a pair $ (\s,\xi)$ with $\s \sim \xi$, it will be useful to interpret each component $(\s_i,\xi_i)$ as  a set of {\em labeled points} in the following sense: if $x \in \xi_i$, then we assign to $x$ the label $\s_i(x^-) = \lim_{t \uparrow x} \s_i(t)$. Note that, if $\xi_i$ is nonempty, then assigning a ``coloring''  $\s_i$ is equivalent to assigning a label to each point of $\xi_i$. It has to be remembered that a label (i.e. a color) is also assigned to the empty configuration of points. It is therefore easy to define the map
\[
\varphi: \S \times \S \times \{-1,1\} \ra \D \times \S,
\]
where, for $(v,w,s)  \in  \S \times \S \times \{-1,1\} $, $\varphi(v,w,s) \in \D \times \S$ is obtained as follows, where $\Pi_{\D} : \D \times \S \ra \D$ and $\Pi_{\S}: \D \times \S \ra \S$ are the canonical projections:
\bi
\item
$\Pi_{\S}\varphi(v,w,s) := v \cup w$.
\item
$\Pi_{\D}\varphi(\emptyset,\emptyset,s) \equiv s$. Whenever $v \cup w \neq \emptyset$, $\Pi_{\D}\varphi(v,w,s)$ is obtained by assigning label $1$ to the points of $v$ and label $-1$ to the points of $w$ (in particular it does not depend on $s$). 
\ei
The function $\varphi$ above can be lifted to a function 
\be{phi}
\Phi: (\S \times \S \times \{-1,1\} )^{\L} \ra S = (\D \times \S)^{\L}
\eeq
 componentwise: for $i \in \L$ and $(\xi,\eta,\rho) \in (\S \times \S \times \{-1,1\} )^{\L} $, we set $(\Phi(\xi,\eta,\rho))_i := \varphi(\xi_i,\eta_i,\rho_i)$. In what follows we set
\[
\hS := (\S \times \S \times \{-1,1\} )^{\L} ,
\]
and 
\[
\hp(d\xi,d\eta,\rho):= \frac{1}{Z} \exp\left[\int_0^{\b} H(\s(\t))d\t \right] \hP(d\xi) \hP(d\eta),
\]
where $Z$ is a normalization factor,  $\s = (\s_i)_{i \in \L}$ is defined by
\be{sigma}
\s_i := \Pi_{\D}\varphi(\xi_i,\eta_i,\rho_i),
\eeq
$\hP := \otimes_{i \in \L} \hat{P}_i$ and $\hat{P}_i$ is the Poisson measure on $\S$ with intensity $\l_i/2$. It is clear that
\be{hp}
\pi = \hp \circ \Phi^{-1}.
\eeq
In the remaining part of this paper we will deal with the probability $\hp$, for which we will establish exponential decay of correlations. By (\ref{hp}) the analogous property for $\pi$ will follow readily.

\subsection{The Markov process} \label{markov}

The aim of this section is to construct a continuous-time Markov process on $\hS$ having $\hp$ as stationary distribution, and show that the spectral gap of its infinitesimal generator has a positive lower bound that does not depend on $\L$.
For this purpose we introduce the following maps from $\hS$ to $\hS$.
\bi
\item {\em Addition maps}.
For $i \in \L$, $x \in [0,\b)$, $\g_{+,1,i}^x (\xi,\eta,\rho)$ is obtained from $(\xi,\eta,\rho)$ by adding the point $x$ to $\xi_i$. Similarly, $\g_{+,2,i}^x (\xi,\eta,\rho)$ is obtained from $(\xi,\eta,\rho)$ by adding the point $x$ to $\eta_i$.
\item {\em Removal maps}.
For $i \in \L$, $x \in [0,\b)$, $\g_{-,1,i}^x (\xi,\eta,\rho)$ is obtained from $(\xi,\eta,\rho)$ by removing the point $x$ to $\xi_i$ (leaving it unchanged if $x \not\in \xi_i$). Similarly, $\g_{-,2,i}^x (\xi,\eta,\rho)$ is obtained from $(\xi,\eta,\rho)$ by removing the point $x$ to $\eta_i$.
\item {\em Spin-flip maps}.
For $i \in \L$, $\g_{0,i} (\xi,\eta,\rho)$ is obtained from $(\xi,\eta,\rho)$ by changing sign to $\rho_i$.
\ei
We now consider the following linear operator, acting on bounded measurable functions from $\hS$ to $\R$.
\begin{multline}
\label{generator}
\gen f(\xi,\eta,\rho) := \sum_{i \in \L} \sum_{k=1,2} \frac{\l_i}{2}  \int_0^{\b} \exp\left[  \nabla_{\g_{+,k,i}^x} \int_0^{\b} H(\s(\t))d\t \right] \nabla_{\g_{+,k,i}^x} f(\xi,\eta,\rho) \\ + \sum_{i \in \L} \sum_{x \in \xi_i}\nabla_{\g_{-,1,i}^x} f(\xi,\eta,\rho) + \sum_{i \in \L} \sum_{x \in \eta_i}  \nabla_{\g_{-,2,i}^x} f(\xi,\eta,\rho) \\ + \sum_{i \in \L}  \exp\left[ \frac{1}{2} \nabla_{\g_{0,i}} \int_0^{\b} H(\s(\t))d\t \right] \nabla_{\g_{0,i}^x}f(\xi,\eta,\rho).
\end{multline}
$\gen$ is the generator of a Markov process that can be described by the following algorithm.
\ben
\item
Let $(\xi,\eta,\rho) \in \hS$ be the state at time $s$. 
\item
Consider three independent families of exponential random variables, $(X_{i,x,k})_{i \in \L, x \in \xi_i, k=1,2}$, $(Y_{i,k})_{i \in \L, k=1,2}$ and $(Z_i)_{i \in \L}$. Within each family, random variables are independent. $X_{i,x,k}$ has mean $\mu^X$, $Y_{i,k}$ has mean $\mu^Y_i$ and $Z_{i}$ has mean $\mu^Z$, where
\[
\frac{1}{\mu^Y_i} = \frac{\l_i}{2} \b e^{\b C}, \ \ \frac{1}{\mu^X} = e^{\b C}, \ \  \frac{1}{\mu^Z} = e^{\b C/2}.
\]
where $C$ is the constant appearing in condition ({\bf BD}).
\item
Let $t$ be the value of the smallest of these random variables. Set $(\xi(s+h), \eta(s+h), \rho(s+h)) \equiv (\xi,\eta,\rho)$ for $0 \leq h < t$.
\bi
\item[a)]
If $t = X_{i,x,k}$, then $(\xi(s+t), \eta(s+t), \rho(s+t)) := \g_{-,k,i}^x (\xi,\eta,\rho)$. 
\item[b)]
If $t = Y_{i,k}$, then $(\tilde{\xi}, \tilde{\eta}, \tilde{\rho}) := \g_{+,k,i}^x (\xi,\eta,\rho)$, where $x$ is sampled from the uniform probability on $[0,\b)$. Let $U$ be a random number uniformly distributed on $[0,1]$. Set
\[
(\xi(s+t), \eta(s+t), \rho(s+t)) := \left\{ \begin{array}{ll} (\tilde{\xi}, \tilde{\eta}, \tilde{\rho}) & \mbox{if } U \leq \frac{ \exp\left[ \nabla_{\g_{+,k,i}^x} \int_0^{\b} H(\s(\t))d\t \right]}{e^{\b C}} \\ (\xi,\eta,\rho) & \mbox{otherwise.} \end{array} \right.
\]
\item[c)]
If $t = Z_{i}$, then $(\tilde{\xi}, \tilde{\eta}, \tilde{\rho}) := \g_{0,i} (\xi,\eta,\rho)$. Set
\[
(\xi(s+t), \eta(s+t), \rho(s+t)) := \left\{ \begin{array}{ll} (\tilde{\xi}, \tilde{\eta}, \tilde{\rho}) & \mbox{if } U \leq \frac{ \exp\left[ \frac{1}{2} \nabla_{\g_{0,i}} \int_0^{\b} H(\s(\t))d\t \right]}{e^{\b C/2}} \\ (\xi,\eta,\rho) & \mbox{otherwise.} \end{array} \right.
\]
\ei
\item
Replace $(\xi, \eta, \rho)$ by $(\xi_{s+t},\eta_{s+t}, \rho_{s+t})$, $s$ by $s+t$, and go to step 1.
\een
It can be shown that the above algorithm provides a rigorous construction of the Markovian dynamics having $\mathcal{L}$ as infinitesimal generator, and that bounded measurable functions form a core for $\gen$. Moreover,  $\gen$ is of the form (\ref{eq:Lintegrale}), once we define 
\[
G := \{\g_{+,k,i}^x, \g_{-,k,i}^x, \g_{0,i}: \, i \in \L, \, x \in [0,\b), \, k = 1,2\},
\]
and $c(\xi,\eta,\rho,d\g)$ by
\begin{multline}
\label{rates}
\int l(\g) c(\xi,\eta,\rho,d\g) := \sum_{i \in \L} \sum_{k=1,2} \frac{\l_i}{2}  \int_0^{\b} \exp\left[ \nabla_{\g_{+,k,i}^x} \int_0^{\b} H(\s(\t))d\t \right] l(\g_{+,k,i}^x) \\ + \sum_{i \in \L} \sum_{x \in \xi_i} l(\g_{-,1,i}^x)  + \sum_{i \in \L} \sum_{x \in \eta_i}  l(\g_{-,2,i}^x)  + \sum_{i \in \L}  \exp\left[ \frac{1}{2} \nabla_{\g_{0,i}} \int_0^{\b} H(\s(\t))d\t \right] l(\g_{0,i}).
\end{multline}
Note that the measurable structure in $G$ is induced by $[0,\b)$: a function $l:G \ra \R$ is measurable if and only if $x \mapsto l(\g_{+,k,i}^x)$ and $x \mapsto l(\g_{-,k,i}^x)$ are measurable. In order to show that property {\bf (Rev)} holds (see (\ref{db})), we first need to observe that, in the sense of {\bf (Rev)}, $\left(\g_{+,k,i}^x \right)^{-1} = \g_{-,k,i}^x$, $\left(\g_{-,k,i}^x \right)^{-1} = \g_{+,k,i}^x$, and $\left(\g_{0,i}\right)^{-1} = \g_{0,i}$.

\bl{reversibility}
The reversibility condition {\bf (Rev)} holds for the generator $\gen$ defined in (\ref{generator}).
\el
\bpr
It will be convenient to define, for $i \in \L$,
\[
G_{i,k}^+ := \{\g_{+,k,i}^x :  x \in [0,\b)\}, \ \  G_{i,k}^-  := \{\g_{-,k,i}^x:  x \in [0,\b)\}, \ \ G_i^0 := \{\g_{0,i}\}.
\]
We will also write $G_i := G_{i,1}^+ \cup G_{i,2}^+ \cup G_{i,1}^- \cup G_{i,2}^- \cup G_i^0$. All these subsets of $G$ are disjoint. Finally, we set
\[
\Hb(\xi,\eta,\rho) := \int_0^{\b} H(\s(\t))d\t,
\]
which makes clearer the following computation, where $\Psi: \hat{S} \times G \ra [0,+\infty)$ is measurable.
\begin{multline*}
\hp\left[\int_{G_{i,1}^+} \Psi(\xi,\eta,\rho,\g) c(\xi,\eta,\rho,d\g)\right] =  \frac{\l_i}{2} \hp\left[\int_0^{\b} dx\exp\left[\nabla_{\g_{+,1,i}^x} \Hb(\xi,\eta,\rho) \right]  \Psi(\xi,\eta,\rho,\g_{+,1,i}^x) \right] \\ =  \frac{\l_i}{2} \frac{1}{Z} \sum_{\rho \in \{\pm 1\}^{\L}} \int \hat{P}_i(d\xi) \hat{P}_i(d\eta) \exp\left[ \Hb(\xi,\eta,\rho) \right] \int_0^{\b} dx \exp\left[  \nabla_{\g_{+,1,i}^x} \Hb(\xi,\eta,\rho) \right]  \Psi(\xi,\eta,\rho,\g_{+,1,i}^x) \\ =  \frac{\l_i}{2} \frac{1}{Z} \sum_{\rho \in \{\pm 1\}^{\L}} \int \hat{P}_i(d\eta) e^{-\l_i/2} \sum_{n=0}^{+\infty} \frac{\l_i^n}{2^n n!} \int_0^{\b} \cdots \int_0^{\b} dx_1 \cdots dx_n \int_0^{\b} dx \\ \exp\left[\Hb(\{x_1,\ldots,x_n,x\}, \eta,\rho) \right] \Psi(\{x_1,\ldots,x_n\}, \eta,\rho,\g_{+,1,i}^x ) \\ = 
\frac{1}{Z} \sum_{\rho \in \{\pm 1\}^{\L}} \int \hat{P}_i(d\eta) e^{-\l_i/2}  \sum_{n=1}^{+\infty} \frac{\l_i^n}{2^n n!} n \int_0^{\b} \cdots \int_0^{\b} dx_1 \cdots dx_n \\ \exp\left[ \Hb(\{x_1,\ldots,x_n\}, \eta,\rho) \right] \Psi(\{x_1,\ldots,x_{n-1}\}, \eta,\rho,\g_{+,1,i}^{x_n} ) \\
= \frac{1}{Z} \sum_{\rho \in \{\pm 1\}^{\L}} \int \hat{P}_i(d\eta) e^{-\l_i/2}  \sum_{n=1}^{+\infty} \frac{\l_i^n}{2^n n!} \int_0^{\b} \cdots \int_0^{\b} dx_1 \cdots dx_n  \\ \exp\left[  \Hb(\{x_1,\ldots,x_n\}, \eta,\rho) \right] \sum_{k=1}^n \Psi(\{x_1,\ldots,x_{n}\} \setminus \{x_k\},\eta,\rho,\g_{+,1,i}^{x_k} ) \\
= \frac{1}{Z} \sum_{\rho \in \{\pm 1\}^{\L}} \int \hat{P}_i(d\eta) e^{-\l_i/2}  \sum_{n=1}^{+\infty} \frac{\l_i^n}{2^n n!} \int_0^{\b} \cdots \int_0^{\b} dx_1 \cdots dx_n  \exp\left[  \Hb(\{x_1,\ldots,x_n\}, \eta,\rho) \right] 
\\ \sum_{k=1}^n 
\Psi\left(\g_{-,1,i}^{x_k}(\{x_1,\ldots,x_{n}\}),\eta,\rho,\left(\g_{-,1,i}^{x_k}\right)^{-1} \right) \\ = \hp\left[\int_{G_{i,1}^-} \Psi(\g(\xi,\eta,\rho),\g^{-1}) c(\xi,\eta,\rho,d\g)\right].
\end{multline*}
Thus, we have shown
\be{rev1}
\hp\left[\int_{G_{i,1}^+} \Psi(\xi,\eta,\rho,\g) c(\xi,\eta,\rho,d\g)\right] = \hp\left[\int_{G_{i,1}^-} \Psi(\g(\xi,\eta,\rho),\g^{-1}) c(\xi,\eta,\rho,d\g)\right].
\eeq
Similarly, one shows
\be{rev2}
\hp\left[\int_{G_{i,1}^-} \Psi(\xi,\eta,\rho,\g) c(\xi,\eta,\rho,d\g)\right] = \hp\left[\int_{G_{i,1}^+} \Psi(\g(\xi,\eta,\rho),\g^{-1}) c(\xi,\eta,\rho,d\g)\right],
\eeq
and, by symmetry, the analogous relations for $G_{i,2}^+, G_{i,2}^-$. Finally, the identity
\be{rev3}
\hp\left[\int_{G_{i}^0} \Psi(\xi,\eta,\rho,\g) c(\xi,\eta,\rho,d\g)\right] = \hp\left[\int_{G_{i}^0} \Psi(\g(\xi,\eta,\rho),\g^{-1}) c(\xi,\eta,\rho,d\g)\right]
\eeq
is simpler to derive, due to the fact that, for $\g \in G_i^0$, we have $\g = \g^{-1}$. Adding up (\ref{rev1}), (\ref{rev2}), the analogous relations for $G_{i,2}^+, G_{i,2}^-$, (\ref{rev3}), and then summing over $i \in \L$, we obtain (\ref{db}).
\epr

\subsection{The spectral gap}

We now apply Theorem \ref{gap} to estimate the spectral gap of $\mathcal{L}$. We recall some notations introduced in the previous section:
\[
G_{i,k}^+ := \{\g_{+,k,i}^x :  x \in [0,\b)\}, \ \  G_{i,k}^-  := \{\g_{-,k,i}^x:  x \in [0,\b)\}, \ \ G_i^0 := \{\g_{0,i}\},
\]
\[
G_i := G_{i,1}^+ \cup G_{i,2}^+ \cup G_{i,1}^- \cup G_{i,2}^- \cup G_i^0
\]
We also set $G_i^+ := G_{i,1}^+ \cup G_{i,2}^+ $, $G_i^- := G_{i,1}^- \cup G_{i,2}^-$.
Consider the following function $r: \hS \times G \times G \ra [0,+\infty)$.
\be{erre}
r(\xi,\eta,\rho,\g,\d):= \left\{
\begin{array}{ll}
	\exp\left[\nabla_{\g} \nabla_{\d} \int_0^{\b} H((\s(\t))d\t \right] & \mbox{for } (\g, \d)  \in G_i^+ \times G_j^+, \, i \neq j \\
	\exp\left[\frac{1}{2}\nabla_{\g} \nabla_{\d} \int_0^{\b} H((\s(\t))d\t \right] & \mbox{for } (\g,\d)  \in \left\{ \begin{array}{l} G_i^+ \times G_j^0 \\ G_i^0 \times G_j^+ \end{array} \right. \, i \neq j \\
	\frac{1}{2} \left[1 + \exp\left[\frac{1}{2}\nabla_{\g} \nabla_{\d} \int_0^{\b} H((\s(\t))d\t \right] \right] & \mbox{for } (\g, \d) \in G_i^0 \times G_j^0, \, i \neq j \\
	1 &  \mbox{for } (\g,\d) \in \left\{ \begin{array}{l} G_i^+ \times G_j^- \\ G_i^- \times G_j^+ \\ G_i^0 \times G_j^- \\ G_i^- \times G_j^0 \\ G_i^- \times G_j^-  \end{array} \right. \, i \neq j \\ 
	0 & \mbox{otherwise}
\end{array} \right.
\eeq
\bp{prop:symmetry}
The function $r$ defined in (\ref{erre}) is $\gen$-symmetric.
\ep
\bpr
Note that, by assumption {\bf (BD)}, $r$ is a bounded function. Thus, condition {\bf (A1)} in Definition  \ref{symmetric} follows if we show that the measure $\hp(d\xi, d\eta, \rho) c(\xi,\eta,\rho,d\g)c(\xi,\eta,\rho,d\d)$ is finite. To see this, note that by (\ref{rates}) and {\bf (BD)}
\[
c(\xi,\eta,\rho,G) \leq e^{\b C} \sum_{i \in \L} \l_i +  \sum_{i \in \L}(|\xi_i| +| \eta_i|) + |\L| e^{\b C/2},
\]
from which
\[
\hp\left[ c(\xi,\eta,\rho,G)c(\xi,\eta,\rho,G) \right] < +\infty
\]
follows. \\
Condition {\bf (A2)} is clearly satisfied by definition of $r$. Condition {\bf (A3)}, is obvious since $\g \circ \d = \d \circ \g$ for every $\g \in G_i, \d \in G_j$, $i \neq j$. We are therefore left to show {\bf (A4)}. Note that, for $F: \hS \times G \times G \ra [0,+\infty)$ measurable,
\begin{multline*}
\int F(\xi,\eta,\rho,\g,\d)r(\xi,\eta,\rho,\g,\d) c(\xi,\eta,\rho,d\g)c(\xi,\eta,\rho,d\d) \hp(d\xi,d\eta,\rho) \\ = \sum_{i \neq j} \hp \left[\int_{G_i \times G_j} F(\xi,\eta,\rho,\g,\d)r(\xi,\eta,\rho,\g,\d) c(\xi,\eta,\rho,d\g)c(\xi,\eta,\rho,d\d) \right].
\end{multline*}
Thus, it is enough to show that for every $i,j \in \L$, $i \neq j$
\begin{multline}
\label{symmetry1}
 \hp \left[\int_{G_i \times G_j} F(\xi,\eta,\rho,\g,\d)r(\xi,\eta,\rho,\g,\d) c(\xi,\eta,\rho,d\g)c(\xi,\eta,\rho,d\d) \right] \\ =  \hp \left[\int_{G_i \times G_j} F(\g(\xi,\eta,\rho),\g^{-1},\d)r(\xi,\eta,\rho,\g,\d) c(\xi,\eta,\rho,d\g)c(\xi,\eta,\rho,d\d) \right].
 \end{multline}
We show (\ref{symmetry1}) by decomposing $G_i$ as $G_i^+ \cup G_i^- \cup G_i^0$. Define, for $A,B \subseteq G$
\begin{eqnarray*}
I(A,B) & := & \hp \left[\int_{A \times B} F(\xi,\eta,\rho,\g,\d)r(\xi,\eta,\rho,\g,\d) c(\xi,\eta,\rho,d\g)c(\xi,\eta,\rho,d\d) \right] \\ 
\hat{I}(A,B) & := & \hp \left[\int_{A \times B} F(\g(\xi,\eta,\rho),\g^{-1},\d)r(\xi,\eta,\rho,\g,\d) c(\xi,\eta,\rho,d\g)c(\xi,\eta,\rho,d\d) \right].
\end{eqnarray*}
Identity (\ref{symmetry1}) follows by adding up the equalities:
\be{s}
I(G^k_i, G^l_j)= \hat I (G^{-k}_i, G^l_j)
\eeq
for $k,l \in \{ -,0,+\}$ (with the obvious meaning for $-k$).
The key fact to prove (\ref{s}) is given in the following two computations, in which $K:\hS \times G \times G \ra [0,+\infty)$ is measurable.
\begin{multline} \label{e4}
\int \hp(d\xi,d\eta,\rho)c(\xi,\eta,\rho,d\g)c(\xi,\eta,\rho,d\d) K(\xi,\eta,\rho,\g,\d) \\ = \int \hp(d\xi,d\eta,\rho)c(\xi,\eta,\rho,d\g)c(\g(\xi,\eta,\rho),d\d) K(\g(\xi,\eta,\rho),\g^{-1},\d),
\end{multline}
where we have applied {\bf (Rev)} to the function $\Psi(\xi,\eta,\rho,\g) := \int c(\eta,\xi,\eta,\rho,d\d) K(\xi,\eta,\rho,\g,\d)$, and
\begin{multline} \label{e5}
\int \hp(d\xi,d\eta,\rho)c(\xi,\eta,\rho,d\g)c(\g(\xi,\eta,\rho),d\d) K(\xi,\eta,\rho,\g,\d) \\= \int \hp(d\xi,d\eta,\rho)c(\xi,\eta,\rho,d\g)c(\xi,\eta,\rho,d\d) K(\g(\xi,\eta,\rho),\g^{-1},\d),
\end{multline}
where {\bf (Rev)}  has been applied to $\Psi(\xi,\eta,\rho,\g) := \int c(\g(\eta,\xi,\eta),\rho,d\d) K(\xi,\eta,\rho,\g,\d)$. \\
We now give some of the proofs of the different cases in (\ref{s}). 

\noindent
{\em Case $(k,l) = (+,+)$}. Observe that, for $\g \in G_i^+$, $\d \in G_j^+$, 
\be{symmetry2}
r(\xi,\eta,\rho,\g,\d) = \frac{dc(\g(\xi,\eta,\rho), \cdot)}{dc(\xi,\eta,\rho,\cdot)}(\d).
\eeq
Thus
\begin{multline*}
I(G_i^+, G_j^+) = \hp\left[\int_{G_i^+ \times G_j^+} c(\xi,\eta,\rho,d\g)c(\xi,\eta,\rho,d\d) r(\xi,\eta,\rho,\g,\d) F(\xi,\eta,\rho,\g,\d) \right] \\ = \hp\left[\int_{G_i^+ \times G_j^+} c(\xi,\eta,\rho,d\g)c(\g(\xi,\eta,\rho),d\d) F(\xi,\eta,\rho,\g,\d) \right]  
= \mbox{ ( by (\ref{e5}) for $K := F {\bf 1}_{G_i^+} {\bf 1}_{G_j^+}$) }\\
= \hp\left[\int_{G_i^- \times G_j^+} c(\xi,\eta,\rho,d\g)c(\xi,\eta,\rho,d\d) F(\g(\xi,\eta,\rho),\g^{-1},\d) \right]  = \hat{I}(G_i^-, G_j^+).
\end{multline*}

\noindent
{\em Cases $(k,l) = (+,0)$ and $(k,l) = (+,-)$}. Again (\ref{symmetry2}) holds, and one proceeds as for the case  $(k,l) = (+,+)$.

\noindent
{\em Case $(k,l) = (0,+)$}. Here (\ref{symmetry2}) does not hold, and is replaced by
\be{symmetry3}
r(\xi,\eta,\rho,\g,\d) = \frac{dc(\d(\xi,\eta,\rho), \cdot)}{dc(\xi,\eta,\rho,\cdot)}(\g).
\eeq
Thus
\begin{multline*}
I(G_i^0, G_j^+) = \hp\left[\int_{G_i^0 \times G_j^+} c(\xi,\eta,\rho,d\g)c(\xi,\eta,\rho,d\d) r(\xi,\eta,\rho,\g,\d) F(\xi,\eta,\rho,\g,\d) \right] \\ = 
\hp\left[\int_{G_i^0 \times G_j^+} c(\d(\xi,\eta,\rho),d\g)c(\xi,\eta,\rho,d\d) F(\xi,\eta,\rho,\g,\d) \right] = \mbox{ ( by (\ref{e5}) with $\d$ and $\g$ exchanged)} \\ = \hp\left[\int_{G_i^0 \times G_j^-} c(\xi,\eta,\rho,d\g)c(\xi,\eta,\rho,d\d) F(\d(\xi,\eta,\rho),\g,\d^{-1}) \right] = \mbox{ ( by (\ref{e4}) and $\g \circ \d = \d \circ \g$)} \\ = \hp\left[\int_{G_i^0 \times G_j^-} c(\xi,\eta,\rho,d\g)c(\g(\xi,\eta,\rho),d\d) F(\g(\d(\xi,\eta,\rho)),\g^{-1},\d^{-1}) \right] = \left( \begin{array}{l}\mbox{ since for $\d \in G_j^-$} \\  \mbox{$ \frac{dc(\g(\xi,\eta,\rho), \cdot)}{dc(\xi,\eta,\rho,\cdot)}(\d) = 1$} \end{array} \right)  \\ =  \hp\left[\int_{G_i^0 \times G_j^-} c(\xi,\eta,\rho,d\g)c(\xi,\eta,\rho,d\d) F(\g(\d(\xi,\eta,\rho)),\g^{-1},\d^{-1}) \right] = \left(\begin{array}{l} \mbox{  by (\ref{e5}) with} \\ \mbox{$\d$ and $\g$ exchanged} \end{array}\right) \\ = \hp\left[\int_{G_i^0 \times G_j^+} c(\d(\xi,\eta,\rho),d\g)c(\xi,\eta,\rho,d\d) F(\g(\xi,\eta,\rho),\g^{-1},\d) \right] = \hat{I}(G_i^0, G_j^+) .
\end{multline*}

\noindent
{\em Case $(k,l) = (0,0)$}. Here we have
\be{symmetry4}
r(\xi,\eta,\rho,\g,\d) = \frac{1}{2} \left[1 +  \frac{dc(\g(\xi,\eta,\rho), \cdot)}{dc(\xi,\eta,\rho,\cdot)}(\d)\right].
\eeq
This gives
\begin{multline*}
I(G_i^0, G_j^0) =\hp\left[\int_{G_i^0 \times G_j^0} c(\xi,\eta,\rho,d\g)c(\xi,\eta,\rho,d\d) r(\xi,\eta,\rho,\g,\d) F(\xi,\eta,\rho,\g,\d) \right] \\ = \frac{1}{2} \hp\left[\int_{G_i^0 \times G_j^0} c(\xi,\eta,\rho,d\g)c(\xi,\eta,\rho,d\d)  F(\xi,\eta,\rho,\g,\d) \right] + \\ \frac{1}{2} \hp\left[\int_{G_i^0 \times G_j^0} c(\xi,\eta,\rho,d\g)c(\g(\xi,\eta,\rho),d\d)  F(\xi,\eta,\rho,\g,\d) \right] = \mbox{ ( by (\ref{e4}) and (\ref{e5})) } \\ = \frac{1}{2} \hp\left[\int_{G_i^0 \times G_j^0} c(\xi,\eta,\rho,d\g)c(\g(\xi,\eta,\rho),d\d)  F(\g(\xi,\eta,\rho),\g^{-1},\d) \right] + \\ \frac{1}{2} \hp\left[\int_{G_i^0 \times G_j^0} c(\xi,\eta,\rho,d\g)c(\xi,\eta,\rho,d\d)  F(\g(\xi,\eta,\rho),\g^{-1},\d) \right] = \hat{I}(G_i^0, G_j^0).
\end{multline*}
The proofs of the remaining cases in (\ref{s}) follow similar arguments, and are omitted
\epr

The above choice of the $\gen$-symmetric function $r$, leads to the following identity:
\begin{multline}
\hat{\pi} \left[ \int c(\xi,\eta,\rho,d\g) c(\xi,\eta,\rho,d\d) (1-r(\xi,\eta,\rho,\g,\d)) \nabla_{\g} f(\xi,\eta,\rho)\nabla_{\d} f(\xi,\eta,\rho) \right] \\ = 
\sum_{i \in \L} \hat{\pi} \left[ \int_{G_i} \int_{G_i} c(\xi,\eta,\rho,d\g) c(\xi,\eta,\rho,d\d) \nabla_{\g} f(\xi,\eta,\rho)\nabla_{\d} f(\xi,\eta,\rho) \right]  \\ + \sum_{i \neq j} \hat{\pi} \left[ \int_{G_i} \int_{G_j} c(\xi,\eta,\rho,d\g) c(\xi,\eta,\rho,d\d) [1-r(\xi,\eta,\rho,\g,\d)]\nabla_{\g} f(\xi,\eta,\rho)\nabla_{\d} f(\xi,\eta,\rho) \right].  \label{term2}
\end{multline}
In order to use Theorem \ref{gap}, we give separate estimates of the two terms in the right hand side of (\ref{term2}).
\bl{lest1}
The following inequality holds:
\begin{multline*}
\left| \sum_{i \neq j} \hat{\pi} \left[ \int_{G_i} \int_{G_j} c(\xi,\eta,\rho,d\g) c(\xi,\eta,\rho,d\d) [1-r(\xi,\eta,\rho,\g,\d)]\nabla_{\g} f(\xi,\eta,\rho)\nabla_{\d} f(\xi,\eta,\rho) \right]\right| \\ \leq 5 \max(1,\b \l) (2R)^d \left(e^{3 \b C} -1\right) \E(f,f),
\end{multline*}
where $\E$ is the Dirichlet for of $\gen$ and $\l$ is the constant appearing in Condition {\bf (BD)}.
\el

\bpr
Fix $i,j \in \L$, with $i \neq j$. Then 
\begin{multline}
\label{offd1}
\hat{\pi} \left[ \int_{G_i} \int_{G_j} c(\xi,\eta,\rho,d\g) c(\xi,\eta,\rho,d\d) [1-r(\xi,\eta,\rho,\g,\d)]\nabla_{\g} f(\xi,\eta,\rho)\nabla_{\d} f(\xi,\eta,\rho) \right]  \\ = \hat{\pi} \left[ \int_{G_i^+} \int_{G_j^+} (\cdots)\right] + \hat{\pi} \left[ \int_{G_i^+} \int_{G_j^0} (\cdots)\right] +  \hat{\pi} \left[ \int_{G_i^0} \int_{G_j^+} (\cdots)\right] +  \hat{\pi} \left[ \int_{G_i^0} \int_{G_j^0} (\cdots)\right] .
\end{multline}
We give separate estimates to the four summands in the right hand side of (\ref{offd1}).
\begin{multline*}
\left|  \hat{\pi} \left[ \int_{G_i^+} \int_{G_j^+}c(\xi,\eta,\rho,d\g) c(\xi,\eta,\rho,d\d) [1-r(\xi,\eta,\rho,\g,\d)]\nabla_{\g} f(\xi,\eta,\rho)\nabla_{\d} f(\xi,\eta,\rho) \right] \right| \\
 = \left|  \frac{\l_i \l_j}{4} \sum_{h,k=1}^2 \hat{\pi} \left[ \int_0^{\b}dx  \int_0^{\b}dy \exp\left[ \nabla_{\g_{+,h,i}^x} \int_0^{\b} H(\s(\t))d\t \right] \exp\left[ \nabla_{\g_{+,k,j}^y} \int_0^{\b} H(\s(\t))d\t \right] \right. \right. \\
\left. \left. \left( 1- \exp\left[ \nabla_{\g_{+,h,i}^x}\nabla_{\g_{+,k,j}^y} \int_0^{\b} H(\s(\t))d\t \right] \right) \nabla_{\g_{+,h,i}^x} f(\xi,\eta,\rho)  \nabla_{\g_{+,k,j}^y} f(\xi,\eta,\rho) \right] \right|
\end{multline*}
By assumptions {\bf (BD)} and {\bf (LOC)},
\[
\exp\left[ \nabla_{\g_{+,h,i}^x} \int_0^{\b} H(\s(\t))d\t \right] \leq e^{\b C}
\]
and 
\[
 \left| 1- \exp\left[ \nabla_{\g_{+,h,i}^x}\nabla_{\g_{+,k,j}^y} \int_0^{\b} H(\s(\t))d\t \right] \right| \leq {\bf 1}_{[0,R]}(|i-j|) \left(e^{2\b C} -1 \right).
 \]
 Thus, using also the inequality $2xy \leq x^2 + y^2$, we obtain
 \begin{multline}
 \label{offd2}
 \left|  \hat{\pi} \left[ \int_{G_i^+} \int_{G_j^+}c(\xi,\eta,\rho,d\g) c(\xi,\eta,\rho,d\d) [1-r(\xi,\eta,\rho,\g,\d)]\nabla_{\g} f(\xi,\eta,\rho)\nabla_{\d} f(\xi,\eta,\rho) \right] \right| \\
 \leq \frac{\b \l}{2}{\bf 1}_{[0,R]}(|i-j|) e^{\b C} \left(e^{2\b C} -1 \right) \hat{\pi} \left[ \frac{\l_i}{2}  \sum_{h=1}^2  \int_0^{\b}dx \exp\left[ \nabla_{\g_{+,h,i}^x} \int_0^{\b} H(\s(\t))d\t \right] \left(\nabla_{\g_{+,h,i}^x} f(\xi,\eta,\rho) \right)^2 \right] \\ + \frac{\b \l}{2}{\bf 1}_{[0,R]}(|i-j|) e^{\b C}\left(e^{2\b C} -1 \right) \hat{\pi} \left[ \frac{\l_j}{2}  \sum_{k=1}^2  \int_0^{\b}dy \exp\left[ \nabla_{\g_{+,k,j}^y} \int_0^{\b} H(\s(\t))d\t \right] \left(\nabla_{\g_{+,k,j}^y} f(\xi,\eta,\rho) \right)^2 \right].
 \end{multline}
 The other terms in (\ref{offd1}) are estimated in the same way, obtaining
 \begin{multline}
 \label{offd3}
  \left|  \hat{\pi} \left[ \int_{G_i^+} \int_{G_j^0}c(\xi,\eta,\rho,d\g) c(\xi,\eta,\rho,d\d) [1-r(\xi,\eta,\rho,\g,\d)]\nabla_{\g} f(\xi,\eta,\rho)\nabla_{\d} f(\xi,\eta,\rho) \right] \right| \\ \leq {\bf 1}_{[0,R]}(|i-j|) \frac{1}{2} e^{\b C/2} \left( e^{\b C} -1 \right) \hat{\pi} \left[ \frac{\l_i}{2}\sum_{h=1}^2  \int_0^{\b}dx \exp\left[ \nabla_{\g_{+,h,i}^x} \int_0^{\b} H(\s(\t))d\t \right] \left(\nabla_{\g_{+,h,i}^x} f(\xi,\eta,\rho) \right)^2 \right] \\ + \l \b {\bf 1}_{[0,R]}(|i-j|) e^{\b C}\left(e^{2\b C} -1 \right) \hat{\pi} \left[\exp\left[ \frac{1}{2}\nabla_{\g_{0,j}} \int_0^{\b} H(\s(\t))d\t \right] \left(\nabla_{\g_{0,j}} f(\xi,\eta,\rho) \right)^2 \right].
  \end{multline}
   \begin{multline}
 \label{offd4}
 \left|  \hat{\pi} \left[ \int_{G_i^0} \int_{G_j^0}c(\xi,\eta,\rho,d\g) c(\xi,\eta,\rho,d\d) [1-r(\xi,\eta,\rho,\g,\d)]\nabla_{\g} f(\xi,\eta,\rho)\nabla_{\d} f(\xi,\eta,\rho) \right] \right| \\ \leq  {\bf 1}_{[0,R]}(|i-j|) \frac{1}{4} e^{\b C/2} \left( e^{\b C} -1 \right) \hat{\pi} \left[\exp\left[ \frac{1}{2}\nabla_{\g_{0,i}} \int_0^{\b} H(\s(\t))d\t \right] \left(\nabla_{\g_{0,i}} f(\xi,\eta,\rho) \right)^2 \right] \\ + {\bf 1}_{[0,R]}(|i-j|) \frac{1}{4} e^{\b C/2} \left( e^{\b C} -1 \right) \hat{\pi} \left[\exp\left[ \frac{1}{2}\nabla_{\g_{0,j}} \int_0^{\b} H(\s(\t))d\t \right] \left(\nabla_{\g_{0,j}} f(\xi,\eta,\rho) \right)^2 \right] ,
 \end{multline}
 while the estimate for $\hat{\pi} \left[ \int_{G_i^0} \int_{G_j^+} (\cdots)\right]$ is obtained from (\ref{offd3}) by exchanging $i$ and $j$. Collecting all these estimates, summing over $i \neq j$ and using the fact that
 \begin{multline*}
 \E(f,f) = \sum_{i \in \L} \hat{\pi} \left[ \frac{\l_i}{2}\sum_{h=1}^2  \int_0^{\b}dx \exp\left[ \nabla_{\g_{+,h,i}^x} \int_0^{\b} H(\s(\t))d\t \right] \left(\nabla_{\g_{+,h,i}^x} f(\xi,\eta,\rho) \right)^2 \right] \\ + \frac{1}{2}  \sum_{i \in \L} \hat{\pi} \left[ \exp\left[ \frac{1}{2}\nabla_{\g_{0,i}} \int_0^{\b} H(\s(\t))d\t \right] \left(\nabla_{\g_{0,i}} f(\xi,\eta,\rho) \right)^2 \right] ,
 \end{multline*}
 the proof is easily completed.
 \epr
 
 Note that the bound in Lemma \ref{lest1} is of order $O(\b)$ as $\b \ra 0$. We now proceed to give a lower estimate to the first term in the right hand side of (\ref{term2}).
\bl{lest2}
The following inequality holds:
\[
\sum_{i \in \L} \hat{\pi} \left[ \int_{G_i} \int_{G_i} c(\s,\xi,d\g) c(\s,\xi,d\d) \nabla_{\g} f(\s,\xi)\nabla_{\d} f(\s,\xi) \right] \geq e^{-2\b C}  \E(f,f).
\]
\el
\bpr
We begin by writing
\begin{multline*}
\sum_{i \in \L} \hat{\pi} \left[ \int_{G_i} \int_{G_i} c(\xi,\eta,\rho,d\g) c(\xi,\eta,\rho,d\d) \nabla_{\g} f(\xi,\eta,\rho)\nabla_{\d} f(\xi,\eta,\rho) \right]  \\ = \sum_{i \in \L} \hat{\pi} \left[ \hat{\pi}^i \left(\int_{G_i} \int_{G_i} c(\xi,\eta,\rho,d\g) c(\xi,\eta,\rho,d\d) \nabla_{\g} f(\xi,\eta,\rho)\nabla_{\d} f\xi,\eta,\rho) \right) \right] =\sum_{i \in \L} \hat{\pi} \left[\hat{\pi}^i \left((\gen_i f)^2 \right) \right],
\end{multline*}
where 
\begin{multline*}
\gen_i f(\xi,\eta,\rho) :=  \int_{G_i}c(\xi,\eta,\rho,d\g)\nabla_{\g} f(\xi,\eta,\rho) = \\ \sum_{k=1,2} \frac{\l_i}{2}  \int_0^{\b} \exp\left[  \nabla_{\g_{+,k,i}^x} \int_0^{\b} H(\s(\t))d\t \right] \nabla_{\g_{+,k,i}^x} f(\xi,\eta,\rho) \\ + \sum_{x \in \xi_i}\nabla_{\g_{-,1,i}^x} f(\xi,\eta,\rho)+ \sum_{x \in \eta_i}  \nabla_{\g_{-,2,i}^x} f(\xi,\eta,\rho) \\ + \exp\left[ \frac{1}{2} \nabla_{\g_{0,i}} \int_0^{\b} H(\s(\t))d\t \right] \nabla_{\g_{0,i}^x}f(\xi,\eta,\rho),
\end{multline*}
and $\hat{\pi}^i$ is obtained by conditioning $\hp$ to $(\xi_j,\eta_j,\rho_j)_{j \neq i}$.

Now, define
\[
\E_i(f,f) := -\hat{\pi}^i \left[ f \gen_i f \right].
\]
We recall that all these $i$-dependent expressions involve conditional expectations, so they depend on the conditioning variables. Suppose now we can find a constant $D>0$ such that the inequality
\be{locgap}
D \hat{\E}_i(f,f)  \leq \hat{\pi}^i \left[ \left( \gen_i f \right)^2 \right]
\eeq
holds for all $i \in \L$, all $f$ measurable and bounded, uniformly on the conditioning. Then, clearly,
\begin{multline} \label{globgap}
\sum_{i \in \L} \hat{\pi} \left[ \int_{G_i} \int_{G_i} c(\xi,\eta,\rho,d\g) c(\xi,\eta,\rho,d\d) \nabla_{\g} f(\xi,\eta,\rho)\nabla_{\d} f(\xi,\eta,\rho) \right]= \sum_{i \in \L} \hat{\pi} \left[\hat{\pi}^i \left((\gen_i f)^2 \right) \right]\\ \geq D \sum_{i \in \L} \hat{\pi}\left[\E_i(f,f) \right] = D \E(f,f).
\end{multline}
We are therefore left to show that (\ref{locgap}) holds with $D = e^{-2\b C}$. We first observe that the conditional measure $\hat{\pi}^i$ can be written as follows:
\be{beta}
\hat{\pi}^i (d\xi_i,d\eta_i,\rho_i) = \frac{\hat{P}_i(d\xi_i) \hat{P}_i(d\eta_i) }{\sum_{s=\pm 1} \exp\left[ \int_0^{\b}(H(\s^{i,s}(\t)) - H(\s(\t)))dt \right] \hat{P}_i(d\xi'_i) \hat{P}_i(d\eta'_i) }
\eeq
where $\s^{i,s}$ denotes the ``coloring'' of the image under $\Phi$ (see (\ref{phi}) and (\ref{sigma})) of the configuration obtained from $(\xi,\eta,\rho)$ replacing $(\xi_i,\eta_i,\rho_i)$ by $(\xi'_i,\eta'_i,s)$, and $\hat{P}_i$ is the Poisson measure on $\S$ with intensity $\frac{\l_i}{2}$. By (\ref{beta}) it is also easy to show that $\gen_i$ is self-adjoint in $L^2(\hp^i)$. Note that, by letting $H \equiv 0$ in (\ref{beta}), we obtain the probability $q_i$ on $\S \times \S \times \{-1,1\}$ given by
\be{betazero}
q_i(dv,dw,s) = \frac{1}{2} \hat{P}_i(dv)\hat{P}_i(dw),
\eeq
and the inequalities
\be{rnbounds}
e^{-\b C} \leq \frac{d\hat{\pi}^i}{dq_i} \leq e^{\b C}
\eeq
hold.
Similarly, letting $H \equiv 0$ in $\genh_i$, we obtain the operator $\A_i$ given by
\begin{multline*}
\A_i f(v,w,s) = \frac{\l_i}{2}\int_0^{\b} dx[ f(v \cup \{x\},w,s) - f(v,w,s)] +  \frac{\l_i}{2}\int_0^{\b} dx[ f(v ,w \cup \{x\} ,s) - f(v,w,s)] \\  
+ \sum_{x \in v} [f(v \setminus \{x\},w,s) - f(v,w,s)] +  \sum_{x \in w} [f(v ,w \setminus \{x\},s) - f(v,w,s)] \\
+  [f(v,w,-s) - f(v,w,s)]  ,
\end{multline*}
which is self-adjoint in $L^2(q_i)$. By comparing $\genh_i$ and $\A_i$, we obtain the bounds
\be{dirbounds}
- e^{-\b C} q_i[f \A_i f] \leq \hat{\E_i}(f,f) \leq -e^{\b C} q_i[f \A_i f].
\eeq

By (\ref{rnbounds}) and (\ref{dirbounds}) we have that the inequality
\be{gap1d}
\g q_i[f;f] \leq - q_i[f \A_i f]
\eeq
implies
\[
\g e^{-2 \b C} \hat{\pi}^i[f;f] \leq - \hat{\E}_i(f,f),
\]
for every $f$ bounded,
which, by Proposition \ref{Gamma2}, implies
(\ref{locgap}) with $D = e^{-2\b C} \g$. Thus, we are left to show that (\ref{gap1d}) holds with $\g = 1$, which is equivalent to
\be{gapa}
\mbox{gap}(\A_i) \geq 1.
\eeq
The Markov process generated by $\A_i$, with initial measure $q_i$, is such that the three components $v_t,w_t,s_t$ are independent; moreover $(v_t)$ and $(w_t)$ are birth and death processes with birth rate $\l_i/2$ and death rate $1$, while $(s_t)$ is a spin-flip process with spin-flip rate $1$. The infinitesimal generators of these processes have all spectral gap equal to $1$; by the {\em tensor property} of spectral gap, it follows that $ \mbox{gap}(\A_i) =1$, which shows (\ref{gapa}).
\epr

Collecting the results in Lemmas \ref{lest1} and \ref{lest2},  we obtain the following lower estimate for the spectral gap of $\mathcal{L}$.
\bt{thmgap}
The following inequality holds:
\[
\mbox{gap}(\mathcal{L}) \geq e^{-2\b C}  -  5 \max(1,\b \l) (2R)^d \left(e^{3 \b C} -1\right) .
\]
\et
\bpr
It is enough to apply the estimates in  Lemmas \ref{lest1} and \ref{lest2} to the two terms in (\ref{term2}), and to use Theorem \ref{gap}.
\epr
Note that the lower bound in Theorem \ref{thmgap} goes to $1$ as $\b \downarrow 0$. In particular it is positive for $\b$ small enough.

\section{Decay of correlations: proof of Theorem \ref{decay}}

We begin by introducing some notations. For $f : \hS \ra \R$ measurable and bounded, and $i \in \L$, we define
\[
M_i(f) :=  \max \left(\sup_{x \in [0,\b)} \sup_{k=1,2} \left\| \nabla_{\g_{+,k,i}^x} f \right\|_{L^{\infty}(\hp)}, \left\| \nabla_{\g_{0,i}} f \right\|_{L^{\infty}(\hp)} \right),
\]
\be{triple}
\L_f := \{ i \in \L : M_i(f) \neq 0\}, \hspace{2cm}  \mnorm{f} := \sum_{i \in \L} M_i(f).
\eeq
Note that the definitions in (\ref{triple}) are extensions of those given in (\ref{norms}) and (\ref{support}), once $\varphi: \O_{\L} \ra \R$ is identified with $f_{\varphi}: \hS \ra \R$ by $f_{\varphi}(\xi, \eta, \rho) := \varphi(\s(0))$, where $\s$ is defined in terms of $(\xi, \eta, \rho)$ in (\ref{sigma}). 
We first prove the so-called {\em finite speed of propagation Lemma}. We set $S_t := e^{t \mathcal{L}}$.
\bl{lemfsp}
Consider two functions $f,g :\hS \ra \R$. Then
\be{finitespeed}
\left|\pi\left[S_t(fg) - S_t f S_t g \right] \right| \leq N \mnorm{f} \mnorm{g} e^{Mt - \e \distanza(\L_f,\L_g)},
\eeq
where 
\[
N = \left(\frac{2}{1-e^{-1/R}}\right)^d , \hspace{2cm} M = 10 \max(1,\l \b) e^{\b C+1} , \hspace{2cm} \e = \frac{1}{2R}.
\]
\el
\bpr
Let $f : \hS \ra \R$ measurable and bounded, and $k \in \L$. Our first aim is to show that there are constants $B,\d>0$, which do not depend on $f$ or $k$, such that
\be{fsp1}
M_k(S_t f) \leq e^{Bt} \sum_{j \in \L} e^{-\d |j-k|} M_j(f).
\eeq
We will also provide explicit values for these constants.
Observe that, for $\g \in G$,
\be{fsp1bis}
\frac{d}{dt} \nabla_{\g}S_t f = \nabla_{\g} \mathcal{L}S_t f  = \mathcal{L} \nabla_{\g}S_t f + \o{[}\nabla_{\g},\mathcal{L}\c{]} S_t f,
\eeq
where $\o{[}\nabla_{\g},\mathcal{L}\c{]} := \nabla_{\g} \mathcal{L} -  \mathcal{L} \nabla_{\g} $ is the commutator between $\nabla_{\g}$ and $\gen$. By a direct computations one obtains the following expressions:
\begin{multline} \label{comm1}
\o{[}\nabla_{\g_{+,h,j}^y},\mathcal{L}\c{]} f(\xi,\eta,\rho) \\ = \sum_{i \in \L} \sum_{k=1,2} \frac{\l_i}{2} \int_0^{\b}dx \nabla_{\g_{+,h,j}^y} \exp\left[ \nabla_{\g_{+,k,i}^x} \int_0^{\b} H(\s(\t)) d\t \right]  \nabla_{\g_{+,k,i}^x} f(\g_{+,h,j}^y(\xi,\eta,\rho)) \\ + \nabla_{\g_{+,h,j}^y} f (\g_{+,h,j}^y(\xi,\eta,\rho)) + \sum_{i \in \L} \nabla_{\g_{+,h,j}^y} \exp\left[ \nabla_{\g_{0,i }} \int_0^{\b} H(\s(\t)) d\t \right] \nabla_{\g_{0,i }} f(\g_{+,h,j}^y(\xi,\eta,\rho)).
\end{multline}
\begin{multline} \label{comm2}
\o{[}\nabla_{\g_{0,j}},\mathcal{L}\c{]}f (\xi,\eta,\rho) \\ = \sum_{i \in \L} \sum_{k=1,2} \frac{\l_i}{2} \int_0^{\b}dx \nabla_{\g_{0,j}} \exp\left[ \nabla_{\g_{+,k,i}^x} \int_0^{\b} H(\s(\t)) d\t \right]  \nabla_{\g_{+,k,i}^x} f(\g_{0,j}(\xi,\eta,\rho)) \\ +  \sum_{i \in \L} \nabla_{\g_{0,j }} \exp\left[ \nabla_{\g_{0,i }} \int_0^{\b} H(\s(\t)) d\t \right] \nabla_{\g_{0,i }} f(\g_{0,j}(\xi,\eta,\rho)).
\end{multline}
Observing that, by {\bf (LOC)},
\begin{multline*}
\nabla_{\g_{+,h,j}^y} \exp\left[ \nabla_{\g_{+,k,i}^x} \int_0^{\b} H(\s(\t)) d\t \right]  =  \nabla_{\g_{+,h,j}^y} \exp\left[ \nabla_{\g_{0,i }} \int_0^{\b} H(\s(\t)) d\t \right] \\  = \nabla_{\g_{0,j}} \exp\left[ \nabla_{\g_{+,k,i}^x} \int_0^{\b} H(\s(\t)) d\t \right] = \nabla_{\g_{0,j }} \exp\left[ \nabla_{\g_{0,i }} \int_0^{\b} H(\s(\t)) d\t \right] = 0
\end{multline*}
for $|i-j| > R$, we obtain the bounds
\be{comm3}
\left\|\o{[}\nabla_{\g_{+,h,j}^y},\mathcal{L}\c{]} f \right\|_{L^{\infty}(\hp)} \leq M_j(f) + 2 \max(1,\l \b) e^{\b C}\sum_{i: |i-j| \leq R} M_i(f),
\eeq
\be{comm4}
\left\|\o{[}\nabla_{\g_{0,j}},\mathcal{L}\c{]} f \right\|_{L^{\infty}(\hp)} \leq 2 \max(1,\l \b) e^{\b C}\sum_{i: |i-j| \leq R} M_i(f).
\eeq
Now observe that, using (\ref{fsp1bis}) for the second equality, for $\g = \g_{+,h,j}^y$ or $\g = \g_{0,j}$, we have
\[
\frac{d}{ds} S_{t-s} \nabla_{\g} S_s f = S_{t-s} \left[ \frac{d}{ds}\nabla_{\g} S_s f - \mathcal{L} \nabla_{\g} S_s f \right] = S_{t-s} \o{[}\nabla_{\g},\mathcal{L}\c{]} S_s f,
\]
or, in integrated form,
\be{fsp7}
 \nabla_{\g} S_t f = S_{t} \nabla_{\g}  f + \int_0^t S_{t-s} \o{[}\nabla_{\g},\mathcal{L}\c{]} S_s f ds
 \eeq
 Using (\ref{comm3}), (\ref{comm3}), (\ref{fsp7}) and the fact that $S_t$ contracts the $\| \cdot \|_{L^{\infty}(\hat{\pi})}$ norm, we obtain
\be{fsp8}
 M_k(S_t f) \leq M_k(f) + \int_0^t \left[ M_k(S_s f) +  4 \max(1,\l \b) e^{\b C}\sum_{i: |i-j| \leq R} M_i(S_s f) \right]ds.
\eeq
Let now $B = (B_{i,j})_{i,j \in \L}$ be the matrix defined by
\[
B_{i,j} := 4 \max(1,\l \b) e^{\b C} {\bf 1}_{\{|i-j| \leq R\}}.
\]
The integral inequality (\ref{fsp8}) implies
\be{fsp9}
M_k(S_t f) \leq e^{t} \sum_{n=0}^{+\infty} \frac{t^n}{n!} \sum_j B^{(n)}_{k,j} M_j(f),
\eeq
where $B^{(n)}_{k,j}$ are the elements of the matrix $B^n$. Set $l := 4 \max(1,\l \b) e^{\b C}$. Thus
\be{fsp10}
B^{(n)}_{k,j} \leq l^n  {\bf 1}_{\{|k-j| \leq n R\}} \leq h^n e^{- \d |k-j|},
\eeq
where, for example
\[
h := e \cdot l, \ \ \d := \frac{1}{R}.
\]
If  we now plug (\ref{fsp10}) in (\ref{fsp9}), we obtain (\ref{fsp1}), with $\d := \frac{1}{R}$ and $B = e l + 1 = 1+ 4 \max(1,\l \b) e^{\b C+1}$.
 
 Having obtained (\ref{fsp1}), we now proceed with the proof of the Lemma. We first observe that, by self-adjointness of $\mathcal{L}$:
\[
2 \E(S_z f, S_z g) = - \pi \left[ S_z f \mathcal{L} S_z g \right] - \pi \left[ S_z g \mathcal{L} S_z f \right] = - \frac{d}{dz}  \pi \left[ S_z f S_z g \right] .
\]
Integrating the previous identity in $[0,t]$ we obtain
\be{propag1}
\pi\left[S_t(fg) - S_t f S_t g \right] = 2 \int_0^t \E(S_z f, S_z g) dz.
\eeq
Since
 \begin{multline*}
 \E(f,g) = \sum_{i \in \L} \hat{\pi} \left[ \frac{\l_i}{2}\sum_{h=1}^2  \int_0^{\b}dx \exp\left[ \nabla_{\g_{+,h,i}^x} \int_0^{\b} H(\s(\t))d\t \right] \nabla_{\g_{+,h,i}^x} f(\xi,\eta,\rho)\nabla_{\g_{+,h,i}^x} g(\xi,\eta,\rho)  \right] \\ + \frac{1}{2}  \sum_{i \in \L} \hat{\pi} \left[ \exp\left[ \frac{1}{2}\nabla_{\g_{0,i}} \int_0^{\b} H(\s(\t))d\t \right] \nabla_{\g_{0,i}} f(\xi,\eta,\rho) \nabla_{\g_{0,i}} g(\xi,\eta,\rho)\right] ,
 \end{multline*}
by (\ref{propag1}) we obtain
\[
\left|\pi\left[S_t(fg) - S_t f S_t g \right] \right| \leq \frac{3}{2} \max(1,\l \b) e^{\b C}  \sum_{i \in \L} \int_0^t M_i(S_s f) M_i(S_s g) ds.
\]
Thus, by (\ref{fsp1}), we have
\begin{multline} \label{propag2}
\left|\pi\left[S_t(fg) - S_t f S_t g \right] \right| \leq  \frac{3}{2} \max(1,\l \b) e^{\b C} \sum_{i \in \L} \sum_{j \in \L_f} e^{-\d |j-i|} M_j(f) \sum_{h \in \L_g} e^{-\d |h-i|} M_h(g)\int_0^t e^{2Bs}ds \\ \leq \frac{3}{4B} \max(1,\l \b) e^{\b C} e^{2Bt} \mnorm{f} \mnorm{g} \sup_{j \in \L_f, h \in \L_g} \sum_{i \in \L} e^{-\d (|j-i | + |h-i|)} \\ \leq \left( \sup_{j,h \in \Z^d} \sum_{i \in \Z^d} e^{-\d (|j-i | + |h-i|)/2 }\right) \frac{3}{4B} \max(1,\l \b) e^{\b C} e^{2 Bt} \mnorm{f} \mnorm{g} e^{-\d \distanza(\L_f,\L_g)/2} \\ \leq C(\d) \frac{3}{4B} \max(1,\l \b) e^{\b C} e^{2Bt} \mnorm{f} \mnorm{g} e^{-\d \distanza(\L_f,\L_g)/2} ,
\end{multline}
where 
\[
C(\d) = \sum_{i \in \Z^d} e^{-\d |i|} \leq \left(\frac{2}{1-e^{-\d}}\right)^d.
\]
Note that, since $B = 1+ 4 \max(1,\l \b) e^{\b C+1}$, then $\frac{3}{4B} \max(1,\l \b) e^{\b C} \leq 1$. Recalling that $\d = 1/R$, by (\ref{propag2}) the inequality (\ref{finitespeed}) follows easily.

\epr

\noindent
{\em Proof of Theorem \ref{decay}}. Let $f,g:\h \ra \R$ be two bounded and measurable functions. By no loss of generality, we assume $\hp[f] = \hp[g] = 0$. By Schwartz inequality
\begin{multline} \label{decay1}
|\hp[f;g] |= |\hp[fg]| = |\hp[S_t (fg)] |= \left|\hp[S_t f S_t g] + \hp\left[S_t(fg) - S_t f S_t g \right]\right| \\  \leq \| S_t f\|_{L^2(\hp)} \| S_t g\|_{L^2(\hp)} + \left|\hp\left[S_t(fg) - S_t f S_t g \right]\right| .
\end{multline}
Set $\g :=  e^{-2\b C}  -  5 \max(1,\b \l) (2R)^d \left(e^{3 \b C} -1\right)$.
By Theorem \ref{thmgap},  $\mbox{gap}(\gen) > \g$. It follows that, $ \| S_t f\|_{L^2(\hp)} \leq  e^{- \g t}\| f\|_{L^2(\hp)}$. Thus, by (\ref{decay1}), and using Lemma \ref{lemfsp}, we obtain
\[
|\hp[f;g] | \leq \| f\|_{L^{\infty}(\hp)} \| g\|_{L^{\infty}(\hp)} e^{-2\g t} +  N \mnorm{f} \mnorm{g} e^{Mt - \e \distanza(\L_f,\L_g)}.
\]
Choosing $t :=  \d \distanza(\L_f,\L_g)/2M$ we obtain, for
\[
\d := \min\left( \frac{1}{2}, \frac{\g}{M} \right)  \ \mbox{and} \ \ \e = \min\left( \frac{1}{2}, \frac{\g}{M} \right)\frac{1}{2R},
\]
\be{decay2}
|\hp[f;g] | \leq \left[ \| f\|_{L^{\infty}(\hp)} \| g\|_{L^{\infty}(\hp)} + N  \mnorm{f} \mnorm{g} \right] \exp\left[- \d  \distanza(\L_f,\L_g) \right].
\eeq
By Theorem \ref{stoch-geom}, using (\ref{decay2}) for $f,g$ only depending of $\s(0)$, we have obtained the inequality stated in Theorem \ref{decay}
\qed

\end{document}